\documentclass[11pt,bezier]{article}
\usepackage{amsmath}
\usepackage{amsfonts,amsthm,amssymb}
\usepackage{amsfonts}
\usepackage{multirow} 
\usepackage{array}
\usepackage{makecell} 

\usepackage{booktabs}
\usepackage[colorlinks, linkcolor=blue, anchorcolor=blue, citecolor=blue]{hyperref}
\usepackage{graphicx}
\usepackage{graphics}
\usepackage{amssymb}
\newtheorem{lem}{Lemma}[section]
\newtheorem{thm}[lem]{Theorem}
\newtheorem{cor}[lem]{Corollary}

\newtheorem{remark}[lem]{Remark}

\theoremstyle{definition}

\begin{document}
\title{Forcibly unicyclic and bicyclic graphic sequences\footnote{The research is supported by National Natural Science Foundation of China (12261086).}}
\author{Peiyi Duan, Yingzhi Tian\footnote{Corresponding author. E-mail: tianyzhxj@163.com (Y. Tian).} \\
{\small College of Mathematics and System Sciences, Xinjiang
University, Urumqi, Xinjiang, 830046, PR China}}

\date{}

\maketitle

\noindent{\bf Abstract}  A sequence $D=(d_1,d_2,\ldots,d_n)$ of non-negative integers  is called  a graphic sequence if there is a simple graph with vertices $v_1,v_2,\ldots,v_n$ such that  the degree of $v_i$ is $d_i$ for $1\leq i\leq n$. Given a graph theoretical property $\mathcal{P}$, a graphic sequence $D$ is forcibly $\mathcal{P}$ graphic if each graph with degree sequence $D$ has property $\mathcal{P}$. A graph is acyclic if it contains no cycles. A connected acyclic graph is just a tree and has $n-1$ edges.
A graph of order $n$ is unicyclic (resp. bicyclic) if it is connected and has $n$ (resp. $n+1$) edges. Bar-Noy, B\"{o}hnlein, Peleg and Rawitz [Discrete Mathematics 346 (2023) 113460] characterized forcibly acyclic and forcibly connected acyclic graphic sequences. In this paper, we aim to characterize forcibly unicyclic and forcibly bicyclic graphic sequences.


\date{}

\maketitle

\section{Introduction}

Let $D$ denote a sequence $(d_1,d_2,\ldots,d_n)$ of non-negative integers. Such a sequence $D$  is the $degree$ $sequence$ of a graph $G$ if the vertices of $G$ can be labeled $v_1,v_2,\ldots,v_n$ such that the degree of $v_i$ is $d_i$ for $1\leq i\leq n$. A sequence is $graphic$ if it is the degree
sequence of some simple graph. Each such graph is a {\it realization} of $D$.  Graphic sequences have many characterizations: early ones include
Havel \cite{Havel}, Erd\"{o}s and Gallai \cite{Erdos}, and Hakimi \cite{Hakimi}.

In order to study what properties of a graph are determined by its degree
sequence, Rao \cite{Rao19812} proposed the concepts of forcibly $\mathcal{P}$ graphic sequence and potentially $\mathcal{P}$ graphic sequence,  where $\mathcal{P}$ is a graph theoretical property.  For a graphic sequence $D$, let $\langle D\rangle$ denote the set of realizations of $D$.  The graphic sequence $D$ is $potentially$ $\mathcal{P}$ $graphic$ if at least one  graph in  $\langle D\rangle$ has property $\mathcal{P}$; the graphic sequence $D$ is $forcibly$ $\mathcal{P}$ $graphic$ if all graphs in $\langle D\rangle$ have property $\mathcal{P}$.

In \cite{Edmonds}, Edmonds characterized potentially $k$-edge-connected graphic sequences, pioneering the research in the field of  potentially graphic sequences. In a series of papers,  Kleitman and Wang \cite{Kleitman, Wang1, Wang2} studied potentially graphic sequences for graph theoretical properties such as $k$-factors, $k$-connectedness and maximal edge-connectedness. Tian, Meng, Lai, and Zhang \cite{Tian} further characterized potentially super edge-connected graphic sequences.

Forcibly planar and forcibly self-complementary graphic sequences were studied by Rao \cite{Rao1978, Rao1981}. Choudum \cite{Choudum} gave some sufficient conditions for a graphic sequence to be forcibly connected. Recently, Liu, Meng and Tian provided some sufficient conditions for a digraphic sequence to be 
forcibly $k$-connected and forcibly $k$-arc-connected \cite{Liu2022};  and investigated a uniform hypergraphic sequence to be forcibly $k$-edge-connected and forcibly super edge-connected  \cite{Liu2023}.

A graph is $acyclic$ if it contains no cycles. Usually, an acyclic graph is also called a $forest$. A connected acyclic graph is a $tree$ and has $n-1$ edges.
A graph of order $n$ is $unicyclic$ (resp. $bicyclic$) if it is connected and has $n$ (resp. $n+1$) edges. Bar-Noy, B\"{o}hnlein, Peleg and Rawitz \cite{Bar-Noy} characterized forcibly acyclic and forcibly connected acyclic graphic sequences. Motivated by those results, we will characterize forcibly unicyclic and forcibly bicyclic graphic sequences.

The structure of this paper is as follows. In Section 2, we will present definitions, terminologies and some results used in the subsequent proofs. Section 3 will characterize the forcibly unicyclic graphic sequences. In the last section, we give a sufficient and necessary condition for a  graphic sequence to be forcibly bicyclic.

\section{Preliminaries}

In this paper,  we assume that all graphs are simple graphs.  For graph-theoretical terminologies and notation not defined here, we follow \cite{Bondy}. Let $G$ be a graph with vertex set $V(G)$ and edge set $E(G)$.
For a vertex $v\in V(G)$, the $neighborhood$ $N_G(v)$ of $v$ in $G$ is the set of vertices adjacent to $v$;   the $closed$ $neighborhood$ $N_G[v]$ of $v$ in $G$ is  $N_G(v)\cup \{v\}$; and the $degree$ $d_G(v)$ of  $v$ in $G$ is $|N_G(v)|$. If $d_G(v)=1$, then $v$ is a $pendant$ $vertex$ of $G$. The $girth$ of $G$, denoted by $girth(G)$, is the length of the shortest cycle in $G$. The $distance$ between two vertices $u$ and $v$ in a connected graph $G$, denoted by $dist_G(u, v)$, is the number of the edges of the shortest path connecting them. The $diameter$ $diam(G)$ of the graph $G$ is the largest distance among all pairs of vertices in $V(G)$.
For the sake of brevity, we usually omit the subscript $G$ in $N_G(v)$, $N_G[v]$, $d_G(v)$ and $dist_G(u, v)$.

Let $V(G)=\{v_1,v_2,\ldots,v_n\}$. Then the $degree$ $sequence$ $D(G)$ of $G$  is defined as the sequence $(d_1,d_2,\ldots,d_n)$, where $d_i=d(v_i)$ for $1\leq i\leq n$. For notational convenience, we often sort this sequence $D(G)=(d_1,d_2,\ldots,d_n)$ in a non-increasing order, that is, $d_1\ge d_2\ge\ldots\ge d_n$. 
Specifically, the sequence is sometimes written in the form $D(G)=(d_1^{c_{1}}, d_2^{c_{2}}, \ldots, d_k^{c_{k}})$, where $c_i$ represents the number of occurrences of $d_i$ in $D(G)$ for $1\leq i\leq k$.

We use $P_k$ for the $path$ with $k+1$ vertices, and $C_n$ and $K_n$ for the $cycle$ and the $complete$ graph with $n$ vertices, respectively. 
The $complete$ $bipartite$ graph is denoted by $K_{s,t}$, where $s$ and $t$ are the number of vertices of its bipartition. 
A $forest$ is an acyclic graph, and a $tree$ is a connected acyclic graph. A tree with diameter at most 2 is a $star$ and a  tree with diameter 3 is a $double$-$star$.
The $sandglass$ is the graph obtained from the disjoint union of two $C_3$s by adding an edge between a vertex in one $C_3$ and a vertex in the other $C_3$. The $generalized$ $sandglass$ $S(r,s,t)$ is the graph obtained from the disjoint union of  $C_r$ and $C_s$ by adding the path $P_t$ connecting a vertex in $C_r$ and a vertex in $C_s$.
The $bowtie$ is the graph obtained from two $C_3$s by identifying a vertex in one $C_3$ and a vertex in the other $C_3$. The $generalized$ $bowtie$ $B(r,s)$ is the graph obtained from $C_r$ and $C_s$ by identifying a vertex in $C_r$ and a vertex in $C_s$.
The $kite$ is the graph obtained from $K_4$ by deleting one edge in $K_4$. The $theta$ $graph$, denoted by $\Theta(r,s,t)$, is the graph obtained  from $P_r$, $P_s$ and $P_t$ by identifying one end vertices of each $P_r$, $P_s$ and $P_t$ first, and then identifying the other end  vertices of each $P_r$, $P_s$ and $P_t$.  
Clearly, $S(3,3,1)$,  $B(3,3)$ and $\Theta(1,2,2)$ are the sandglass graph, the bowtie graph and the kite graph, respectively.

In the proofs of this paper, we often use the $switching$ $operation$. Switching $\{uv, xy\}$ with $\{ux,vy\}$ in $G$ is defined as  deleting the two edges $uv$, $xy$ from $G$ and adding two edges $ux$, $vy$ to $G$, where $ux$ and $vy$ do not appear in the graph $G$ originally. Of course, we can switch more than two edges.
Note that, by switching operation, we obtained a graph with the same degree sequence as $G$. 

\begin{lem}\label{lem21}
If a connected graph $G$ has two pendant vertices with distance greater than $3$, then there exists a new graph $G'$ with the same degree sequence as $G$ such that it has exactly two components, one of which is an isolated edge.
 \end{lem}

\noindent{\bf Proof.} Let $u$ and $v$ be two pendant vertices with distance greater than $3$ in $G$. Assume $u$ is adjacent to $u'$ and $v$ is adjacent to $v'$. Then, by switching $\{uu', vv'\}$ with $\{uv,u'v'\}$ in $G$, we obtain a new graph $G'$ with the same degree sequence as $G$ such that it has exactly two components, one of which is an isolated edge. 
$\hfill\Box$

\begin{cor}\label{cor21}
Let $D$ be a forcibly unicyclic or a forcibly bicyclic graphic sequence, and let $G\in \langle D \rangle$. Then any pair of pendant vertices $u$ and $v$ in $G$ satisfy $dist(u,v)\leq3$.
\end{cor}

In \cite{Bar-Noy}, the authors characterized the forcibly tree graphic sequences.

\begin{thm}(\cite{Bar-Noy})\label{th21}
Let $T$ be a tree. Then $D(T)$ is a forcibly tree graphic sequence if and only if $T$ is a star or a double-star.
\end{thm}

Let $D$ be a forcibly unicyclic or a forcibly bicyclic graphic sequence, and let $G\in \langle D \rangle$. By deleting all edges in the cycles of $G$, we obtain a forest with the vertex set $V(G)$ and denote by $Tx$ 
as the component containing $x$, where $x\in V(G)$. Since $D$ is a forcibly unicyclic or a forcibly bicyclic graphic sequence,  we have $D(Tx)$ is a forcibly tree graphic sequence and $Tx$ is a star or a double-star for any $x\in V(G)$. For otherwise, we can find a disconnected realization of $D$. So the following corollary is obtained.

\begin{cor}\label{cor22}
Let $D$ be a forcibly unicyclic or a forcibly bicyclic graphic sequence, and let $G\in \langle D \rangle$. 
If we delete all edges in the cycles of $G$, then each component of the obtained graph is isomorphic to  a star or a double-star.
\end{cor}

\section{Forcibly unicyclic graphic sequences}

\begin{lem}\label{lem31} 
Let $D=(d_1,d_2,\ldots,d_n)$ be a forcibly unicyclic graphic sequence.  

$(1)$ If $n=3$, then $D=(2,2,2)$. 
 
$(2)$ If $n=4$, then $D=(2,2,2,2)$, or $D=(3,2,2,1)$. 

$(3)$ If $n=5$, then $D=(2,2,2,2,2)$, or $D=(3,2,2,2,1)$, or $D=(3,3,2,1,1)$, or $D=(4,2,2,1,1)$. 
 \end{lem}

\noindent{\bf Proof.} Let $G\in \langle D \rangle$.   

(1) If $n=3$, then $G$ is isomorphic to the 3-cycle $C_3$. So $D=(2,2,2)$.  

(2) If $n=4$, then $G$ is isomorphic to the 4-cycle $C_4$, or the graph obtained from $C_3$ by adding a pendant vertex. So $D=(2,2,2,2)$, or $D=(3,2,2,1)$.   

(3) If $n=5$, then $G$ is isomorphic to the 5-cycle $C_5$, or the graph obtained from $C_4$ by adding a pendant  vertex, or the graph obtained from $C_3$ by adding two pendant  vertices (the two pendant  vertices can be adjacent to two different vertices or the same vertex of  $C_3$), or the graph obtained from the disjoint union of $C_3$ and $K_2$ by adding an edge between $C_3$ and $K_2$. So $D=(2,2,2,2,2)$, or $D=(3,2,2,2,1)$, or $D=(3,3,2,1,1)$, or $D=(4,2,2,1,1)$. 
  $\hfill\Box$

\begin{lem}\label{lem32}
If $D=(d_1,d_2,\ldots,d_n)$ is a forcibly unicyclic graphic sequence, then, for any graph $G\in \langle D \rangle$, we have $girth(G)\leq 5$.
 \end{lem}

\noindent{\bf Proof.} On the contrary, assume $girth(G)\ge 6$.
Let $C_k = x_{1}x_{2}\cdots x_{k - 1}x_{k}x_{1}$ be a cycle of $G$,  where $k\ge6$. Then the graph $G'$ obtained from $G$ by switching  $\{x_1x_2,x_4x_5\}$ with  $\{x_1x_5, x_2x_4\}$ has the degree sequence $D$, but it is not connected, a contradiction.  $\hfill\Box$

\begin{lem}\label{lem33}
Let $D=(d_1,d_2,\ldots,d_n)$ be a forcibly unicyclic graphic sequence. If $n\geq6$, then there is a graph $G'\in \langle D \rangle$ such that $girth(G')=3$.
 \end{lem}

\noindent{\bf Proof.} Let $G\in \langle D \rangle$. By Lemma \ref{lem32},  we have $girth(G)\leq 5$. If $girth(G)=3$, then we are done. So assume $4\leq girth(G)\leq 5$. 
Let $C_k = x_{1}x_{2}\cdots x_{k - 1}x_{k}x_{1}$ be a cycle of $G$,  where $4\leq k\leq 5$.  Since $n\geq6$, there is a vertex $y\in V(G)\setminus V(C_k)$ such that $y$ is adjacent to some vertex in $V(C_k)$, say $x_1$. Then the graph $G'$ obtained from $G$ by switching $\{x_1y,x_3x_4\}$ with $\{x_1x_3, x_4y\}$ has the degree sequence $D$ and $girth(G')=3$.  $\hfill\Box$

\begin{lem}\label{lem34}
Let $D=(d_1,d_2,\ldots,d_n)$ be a forcibly unicyclic graphic sequence. If $n\geq6$, then $d_6=1$.
 \end{lem}

\noindent{\bf Proof.}
By Lemma \ref{lem33}, there is a graph $G\in \langle D \rangle$ such that $girth(G)=3$. Let $C_3=x_1x_2x_3x_1$ be a 3-cycle in $G$.
On the contrary, assume $d_6\ge 2$. Then there are three vertices  $y_1, y_2, y_3\in V(G)\setminus \{x_1,x_2,x_3\}$ such that their degrees are all greater than 1.

If $y_1, y_2, y_3$ belong to $Tx_j-x_j$ for some $1\leq j\leq 3$, then  $Tx_j$ is not a star or double-star, which contradicts to Corollary \ref{cor22}. So there are two vertices in $\{x_1,x_2,x_3\}$, say $x_1$ and $x_2$, such that $V(Tx_1-x_1)\cap \{y_1, y_2, y_3\}\neq \emptyset$  and $V(Tx_2-x_2)\cap \{y_1, y_2, y_3\}\neq \emptyset$. Then the distance between a pendant vertex in $Tx_1-x_1$ and a pendant vertex in  $Tx_2-x_2$ is at least 5, which contradicts to Corollary \ref{cor21}. 
$\hfill\Box$

\begin{lem}\label{lem35}
Let $D=(d_1,d_2,\ldots,d_n)$ be a forcibly unicyclic graphic sequence. If $n\geq6$, then $d_4\le 2$.
 \end{lem}

\noindent{\bf Proof.} By Lemma \ref{lem33}, there is a graph $G\in \langle D \rangle$ such that $girth(G)=3$. Let $C_3=x_1x_2x_3x_1$ be a 3-cycle in $G$.
On the contrary, assume $d_4\ge 3$. By Corollary \ref{cor22}, $Tx_i$ contains at most two vertices with degrees greater than two for $i=1,2,3$. Since $d_4\ge 3$, one of $Tx_i$, say $Tx_1$, contains exactly two vertices with degrees greater than 2, and another $Tx_j$, say $Tx_2$, contains at least one vertex with degrees greater than 2. Then we can find a pendant vertex in $Tx_1-x_1$ and a pendant vertex in  $Tx_2-x_2$ such that their distance is at least 4, which contradicts to Corollary \ref{cor21}. 
$\hfill\Box$

\begin{thm}\label{th31}
Let $D=(d_1,d_2,\ldots,d_n)$ be a non-negative and non-increasing integer sequence. Then
$D$ is a forcibly unicyclic graphic sequence if and only if one of the following holds:

$(1)$ $D=(2^5)$ or $D=(3,2^4,1)$;

$(2)$ $D=(n-2,2^3,1^{n-4})$, where $n\ge 4$;

$(3)$ $D=(r,s,t,1^{n-3})$, where $n\ge 3$, $r\geq s\geq t\geq2$ and $r+s+t=n+3$.
\end{thm}

\noindent{\bf Proof.} First, we show the if part. Clearly, all of the sequences in (1), (2) and (3) are graphic sequences and satisfy $\Sigma_{i=1}^n d_i=2n$. It is sufficient to prove that $G$ is connected for any graph $G\in \langle D \rangle$.

$(1)$ Apparently, the sequence $D=(2^5)$ has only one realization, which is $C_5$. Let $G\in \langle (3,2^4,1)\rangle$. Assume $d(u)=3$, $N(u)=\{u_1,u_2,u_3\}$ and $V(G)\setminus N[u]=\{x,y\}$.  Since there is only one vertex with degree 1 in $G$, one of $x$ and $y$,  say $x$, has degree at least 2. Then $x$ is adjacent to some vertex in  $\{u_1,u_2,u_3\}$ and $y$ is adjacent to some vertex in $\{u_1,u_2,u_3, x\}$ 
Therefore, the graph $G$ is connected.

$(2)$  Let $G\in \langle (n-2,2^3,1^{n-4})\rangle$. For $n=4$, $D=(2^4)$ has an unique realization $C_4$. Assume $n>4$. Let $d(u)=n-2$, $N(u)=\{u_1,\cdots,u_{n-2}\}$ and $V(G)\setminus N[u]=\{x\}$. Since $d(x)=1$ or $2$,  we know that $x$  is adjacent to some vertex in $N(u)$. Consequently, $G$ is connected.

$(3)$ Let $G\in \langle (r,s,t,1^{n-3})\rangle$ and $d(x)=r, d(y)=s$, $d(z)=t$. Since all vertices in
$V(G)\setminus\{x,y,z\}$ have degree 1, it follows that $x$, $y$ and $z$ must be adjacent to at least $r-2$, $s-2$ and $t-2$ vertices with degree one, respectively. Since $r+s+t=n+3$ and $G$ has exactly $n$
vertices, we conclude that $x$, $y$ and $z$ are adjacent to exactly $r-2$, $s-2$ and $t-2$ vertices with degree 1, respectively. Therefore, $G$ is the graph obtained from the 3-cycle $xyz$ by adding $r-2$, $s-2$ and $t-2$ pendant vertices to $x$, $y$ and $z$, respectively. Thus $G$ is connected.

Next, we prove the only if part. Let $D$ be a forcibly unicyclic graphic sequence. Then $\Sigma_{i=1}^n d_i=2n$ and $d_1\geq d_2\geq d_3\geq2$.
If $n\leq5$, then, by Lemma \ref{lem31},  $D=(2,2,2)$, or $D=(2,2,2,2)$, or $D=(3,2,2,1)$, or $D=(2,2,2,2,2)$, or $D=(3,2,2,2,1)$, or $D=(3,3,2,1,1)$, or $D=(4,2,2,1,1)$.  The result holds. So we assume $n\geq 6$ in the following. By Lemmas \ref{lem34} and \ref{lem35}, we have $d_6=1$ and $d_4\le2$.
Set $d_1=r$, $d_2=s$ and $d_3=t$.   We consider two cases as follows.

\noindent{\bf Case 1.}  $d_4=1$.

Since $\Sigma_{i=1}^n d_i=2n$, we have $r+s+t=2n-(n-3)=n+3$.  Thus $D=(r,s,t,1^{n-3})$,  where $n\ge 3$, $r\geq s\geq t\geq2$ and $r+s+t=n+3$.

\noindent{\bf Case 2.} $d_4=2$.

By $d_4 = 2$,  we  have  $d_5=1$ or 2.

\noindent{\bf Subcase 2.1.}  $d_5=1$.

Since $\Sigma_{i=1}^n d_i=2n$, we have $r+s+t=2n-2-(n-4)=n+2$.  Thus $D=(r,s,t,2,1^{n-4})$, where $r+s+t=n+2$. 
If $s=2$, then $ t=2$ and $D=(n-2,2^3,1^{n-4})$. 
Suppose $s\ge3$. Then the graph  $G_1$ obtained from  the disjoint union of $\Theta(1,2,2)$ and $K_2$ by adding $r-3$, $s-3$ and $t-2$ pendant vertices to the two vertices with degree 3 and one vertex with degree 2 in $\Theta(1,2,2)$, respectively, has degree sequence $(r,s,t,2,1^{n-4})$. But $G_1$ is not connected, which contradicts that $D$ is an unicyclic graphic sequence.

\noindent{\bf Subcase 2.2.} $d_5=2$.

Since $\Sigma_{i=1}^n d_i=2n$, we have $r+s+t=2n-2-2-(n-5)=n+1$.  Thus $D=(r,s,t,2^2,1^{n-5})$, where $r+s+t=n+1$. If $s=2$, then $r=n-3$  and $D=(n-3,,2^4,1^{n-5})$. Furthermore, if $n\geq7$, then the graph  $G_2$ obtained from  the disjoint union of the bowtie $B(3,3)$ and $K_2$ by adding $n-7$ pendant vertices to the vertex with degree 4 in the  $B(3,3)$ has degree sequence $D=(n-3,,2^4,1^{n-5})$. But $G_2$ is not connected, which contradicts that $D$ is an unicyclic graphic sequence.  So $n=6$ and $D=(3,2^4,1)$. Now we assume $s\geq3$. Then the graph  $G_3$ obtained from  the disjoint union of $\Theta(1,2,3)$ and $K_2$ by adding $r-3$, $s-3$ and $t-2$ pendant vertices to the two vertices with degree 3 and one vertex with degree 2 in $\Theta(1,2,3)$, respectively, has degree sequence $(r,s,t,2^2,1^{n-5})$. But $G_3$ is not connected, which also contradicts that $D$ is an unicyclic graphic sequence. 
$\hfill\Box$

\section{Forcibly bicyclic graphic sequences}

For a bicyclic graph,  it either contains a generalized sandglass graph, or a generalized bowtie graph, or a theta graph.

\begin{lem}\label{lem41} Let $D=(d_1,d_2,\ldots,d_n)$ be  a forcibly bicyclic graphic sequence. If $D$ has a realization $G$ containing  a generalized sandglass graph $S(r,s,t)$, then there exists another realization $G'\in \langle D \rangle$ containing a theta graph $\Theta(3, r+s-3, t)$.
\end{lem}

\noindent{\bf Proof.} 
Let $C_r=x_1x_2\cdots x_rx_1$ and $C_s=y_1y_2\cdots y_sy_1$ be the two disjoint cycles in $G$  and $P_t=z_1\cdots z_{t+1}$ be the unique path connecting $C_r$ and $C_s$ in $G$, where $z_1=x_1$ and $z_{t+1}=y_1$. Then the graph $G'$ obtained from $G$ by switching $\{x_2x_3, y_2y_3\}$ with $\{x_2y_2, x_3y_3\}$  has degree sequence $D$ and contains the theta graph $\Theta(3, r+s-3, t)$.
$\hfill\Box$

\begin{lem}\label{lem42}
Let $D=(d_1,d_2,\ldots,d_n)$ be a forcibly bicyclic graphic sequence. If $D$ has a realization $G$ containing a generalized bowtie graph $B(r,s)$, then there exists another realization $G'\in \langle D \rangle$ containing a generalized bowtie graph $B(3,r+s-3)$.
 \end{lem}
 
\noindent{\bf Proof.}
Let $C_r=x_1x_2\cdots x_rx_1$ and $C_s=y_1y_2\cdots y_sy_1$ be the two cycles in $G$, where $x_1=y_1$. By switching  $\{x_2x_3, y_2y_3\}$ with $\{x_2y_2, y_3x_3\}$, the obtained graph  $G'$  has degree sequence $D$ and contains a generalized  bowtie graph $B(3,r+s-3)$.
$\hfill\Box$

\begin{lem}\label{lem43} 
Let $D=(d_1,d_2,\ldots,d_n)$ be a forcibly bicyclic graphic sequence. If $D$ has a realization $G$ containing a theta graph $\Theta(r,s,t)$, then there exists another realization $G'\in \langle D \rangle$ containing a theta graph $\Theta(1,2,r+s+t-3)$.
 \end{lem}

\noindent{\bf Proof.}
Let $P_r=x_1x_2x_3\cdots x_rx_{r+1}$, $P_s=y_1y_2y_3\cdots y_sy_{s+1}$ and $P_t=z_1z_2z_3\cdots z_tz_{t+1}$ be the three paths in $\Theta(r,s,t)$, where $x_1=y_1=z_1$ and $x_{r+1}=y_{s+1}=z_{t+1}$.
If min$\{r,s,t\}\geq2$, then by switching $\{x_1x_2, y_sy_{s+1}\}$ with $\{x_1y_{s+1}, x_2 y_s\}$, we obtain a graph with degree sequence $D$ and containing $\Theta(1,r+s-1,t)$. So we can assume $r=1$. If min$\{s,t\}\geq3$, then by switching $\{y_2y_3, z_tz_{t+1}\}$ with $\{y_2z_{t+1}, y_3 z_t\}$, we obtain a graph $G'$ with degree sequence $D$ and containing $\Theta(1,2,r+s+t-3)$.
$\hfill\Box$

\begin{remark}\label{rem41}
From the three lemmas above, we can find a realization $G$ in a forcibly bicyclic sequence $D$
such that $G$ either contains a generalized bowtie graph $B(3,s)$, or contains  a theta graph $\Theta(1, 2,t)$. 
\end{remark}

\begin{lem}\label{lem45}
Let $D=(d_1,d_2,\ldots,d_n)$ be a forcibly bicyclic graphic sequence. If $D$ has a realization $G$ containing a generalized bowtie graph $B(3,s)$, then $3\leq s\leq5$.
 \end{lem}

\noindent{\bf Proof.} 
Let $C_3=x_1x_2x_3x_1$ and $C_s=y_1y_2\cdots y_sy_1$ be the two cycles of $B(3,s)$, where $x_1=y_1$. 
If $s\geq6$, then by switching $\{y_2y_3, y_5y_6\}$ with $\{y_2y_6, y_3y_5\}$, we obtain a disconnected realization of $D$, which contradicts that $D$ is forcibly bicyclic.
$\hfill\Box$

\begin{lem}\label{lem46}
Let $D=(d_1,d_2,\ldots,d_n)$ be a forcibly bicyclic graphic sequence. If $D$ has a realization $G$ containing a theta graph $\Theta(1,2,t)$, then $2\leq t\leq4$.
 \end{lem}
 
\noindent{\bf Proof.}
Let $P_r=x_1x_2$, $P_s=y_1y_2y_3$ and $P_t=z_1z_2z_3\cdots z_tz_{t+1}$ be the three paths in $\Theta(1,2,t)$, where $x_1=y_1=z_1$ and $x_2=y_3=z_{t+1}$.
If $t\geq5$, then by switching $\{z_1z_2, z_4z_5\}$ with $\{z_1z_5, z_2z_4\}$, we obtain a disconnected realization of $D$, which contradicts that $D$ is forcibly bicyclic.
$\hfill\Box$

\begin{lem}\label{lem47} 
Let $D=(d_1,d_2,\ldots,d_n)$ be a forcibly bicyclic graphic sequence.  

(1) If $n=4$, then $D=(3,3,2,2)$. 
 
(2) If $n=5$, then $D=(3,3,2,2,2)$, or $D=(3,3,3,2,1)$, or $D=(4,2,2,2,2)$, or $D=(4,3,2,2,1)$. 

(3) If $n=6$, then $D=(3,3,2,2,2,2)$, or $D=(3,3,3,2,2,1)$, or $D=(4,2,2,2,2,2)$, or $D=(4,3,2,2,2,1)$, or $D=(4,3,3,2,1,1)$, or $D=(4,4,2,2,1,1)$, or $D=(5,2,2,2,2,1)$, or $D=(5,3,2,2,1,1)$. 
 \end{lem}

\noindent{\bf Proof.} By remark \ref{rem41}, we can choose  $G\in \langle D \rangle$ such that $G$ contains a generalized bowtie $B(3,s)$ or a theta graph $\Theta(1,2,t)$, where $3\leq s\leq5$ and $2\leq t\leq4$.  

(1) If $n=4$, then $G$ is isomorphic to the theta graph $\Theta(1,2,2)$. So $D=(3,3,2,2)$.  

(2) If $n=5$, then $G$ is isomorphic to  the generalized bowtie $B(3,3)$, or the theta graph $\Theta(1,2,3)$, or the graph obtained from the theta graph $\Theta(1,2,2)$ by adding a pendant vertex. So $D=(3,3,2,2,2)$, or $D=(3,3,3,2,1)$, or $D=(4,2,2,2,2)$, or $D=(4,3,2,2,1)$.

(3) If $n=6$, we will consider two cases as follows.

\noindent{\bf Case 1.}  $G$ contains a generalized bowtie $B(3,s)$, where $3\leq s\leq5$.

If $G$ contains the  generalized bowtie $B(3,4)$, then $G$ is isomorphic to  $B(3,4)$ and  $D=(4,2,2,2,2,2)$.

If $G$ contains the  generalized bowtie $B(3,3)$, then $G$ is obtained from $B(3,3)$ by adding a pendant vertex. So $D=(4,3,2,2,2,1)$, or $D=(5,2,2,2,2,1)$. 

\noindent{\bf Case 2.}  $G$ contains  a theta graph $\Theta(1,2,t)$, where $2\leq t\leq4$.

If $G$ contains the theta graph $\Theta(1,2,4)$, then $G$ is isomorphic to  $\Theta(1,2,4)$  and $D=(3,3,2,2,2,2)$.

If $G$ contains the theta graph $\Theta(1,2,3)$, then $G$ is obtained from $\Theta(1,2,3)$ by adding a pendant vertex. So $D=(3,3,3,2,2,1)$, or $D=(4,3,2,2,2,1)$. 

If $G$ contains the theta graph $\Theta(1,2,2)$, then there are exactly two vertices $x$ and $y$ in $V(G)\setminus V(\Theta(1,2,2))$. When $x$ and $y$ are not adjacent in $G$, then $G$ is obtained from $\Theta(1,2,2)$ by adding two pendant vertices $x$ and $y$.  If $x$ and $y$ are adjacent to different vertices in $\Theta(1,2,2)$, then $D=(4,3,3,2,1,1)$ or $D=(4,4,2,2,1,1)$.   If $x$ and $y$ are adjacent to the same vertex in $\Theta(1,2,2)$, then $D=(4,3,3,2,1,1)$ or  $D=(5,3,2,2,1,1)$.  When $x$ and $y$ are adjacent in $G$, then $G$ is obtained from the union of  $\Theta(1,2,2)$ and $G[\{x,y\}]$ by adding an edge between them.  So  $D=(3,3,3,2,2,1)$ or $D=(4,3,2,2,2,1)$.
$\hfill\Box$

\begin{lem}\label{lem48}
Let $D=(d_1,d_2,\ldots,d_n)$ be a forcibly unicyclic graphic sequence and $D\neq(4,2^6)$. If $n\geq7$, then there is a graph $G'\in\langle D \rangle$ such that $G'$ contains either the bowtie $B(3,3)$, or the kite $\Theta(1,2,2)$ as subgraph.
\end{lem}

\noindent{\bf Proof.}
By Remark \ref{rem41}, Lemmas \ref{lem45} and \ref{lem46}, we can choose  $G\in \langle D \rangle$ such that $G$ contains a generalized bowtie $B(3,s)$ or a theta graph $\Theta(1,2,t)$, where $3\leq s\leq5$ and $2\leq t\leq4$.

\noindent{\bf Case 1.}  $G$ contains a generalized bowtie $B(3,s)$, where $3\leq s\leq5$.

If $s=3$, then we are done.

If $s=4$, then let $ux_1x_2u$ and $uy_1y_2y_3u$ be the two cycles of $B(3,4)$. Since $n\geq7$, there is a vertex $z\in V(G)\setminus V(B(3,4))$ such that $z$ is adjacent to one vertex in $B(3,4)$. By symmetry, we only consider that  $z$ is adjacent to $u$, or $x_1$, or $y_1$ or $y_2$.  

Suppose $z$ is adjacent to $u$. By switching $\{y_1y_2,y_2y_3,uz\}$ with $\{y_1y_3, y_2u, y_2z\}$, we obtain a graph $G_1\in\langle D \rangle$, which contains the bowtie $B(3,3)$ as subgraph.

Suppose $z$ is adjacent to $x_1$. By switching $\{x_1z, y_1y_2,y_2y_3\}$ with $\{xy_2, y_1y_3, y_2z\}$, we obtain a graph $G_2\in\langle D \rangle$, which contains the bowtie $B(3,3)$ as subgraph.

Suppose $z$ is adjacent to $y_1$. By switching $\{y_1z,y_2y_3\}$ with $\{y_1y_3, y_2z\}$, we obtain a graph $G_3\in\langle D \rangle$, which contains the bowtie $B(3,3)$ as subgraph.

Suppose $z$ is adjacent to $y_2$. By switching $\{y_1u,y_2z\}$ with $\{y_1z, y_2u\}$, we obtain a graph $G_4\in\langle D \rangle$, which contains the bowtie $B(3,3)$ as subgraph.

If $s=5$, then let $ux_1x_2u$ and $uy_1y_2y_3y_4u$ be the two cycles of $B(3,5)$. Since $n\geq7$ and  $D\neq(4,2^6)$, there is a vertex $z\in V(G)\setminus V(B(3,5))$ such that $z$ is adjacent to one vertex in $B(3,5)$. By symmetry, we only consider that  $z$ is adjacent to $u$, or $x_1$, or $y_1$ or $y_2$.  

Suppose $z$ is adjacent to $u$. By switching $\{y_2y_3,uz\}$ with $\{y_2u, y_3z\}$, we obtain a graph $G_5\in\langle D \rangle$, which contains the bowtie $B(3,3)$ as subgraph.

Suppose $z$ is adjacent to $x_1$. By switching $\{x_1z, y_1y_2,y_3y_4\}$ with $\{x_1y_2, y_1y_4, y_3z\}$, we obtain a graph $G_6\in\langle D \rangle$, which contains the bowtie $B(3,3)$ as subgraph.

Suppose $z$ is adjacent to $y_1$. By switching $\{y_1z,y_3y_4\}$ with $\{y_1y_4, y_3z\}$, we obtain a graph $G_7\in\langle D \rangle$, which contains the bowtie $B(3,3)$ as subgraph.

Suppose $z$ is adjacent to $y_2$. By switching $\{y_2z,y_4u\}$ with $\{y_2u, y_3z\}$, we obtain a graph $G_8\in\langle D \rangle$, which contains the bowtie $B(3,3)$ as subgraph.

\noindent{\bf Case 2.}  $G$ contains  a theta graph $\Theta(1,2,t)$, where $2\leq t\leq4$.

If $t=2$, then we are done.

If $t=3$, then let $uv$, $ux_1v$ and $uy_1y_2v$ be the three paths of $\Theta(1,2,3)$. Since $n\geq7$, there is a vertex $z\in V(G)\setminus V(\Theta(1,2,3))$ such that $z$ is adjacent to one vertex in $\Theta(1,2,3)$. By symmetry, we only consider that  $z$ is adjacent to $u$, or $x_1$, or $y_1$.  

Suppose $z$ is adjacent to $u$. By switching $\{y_1y_2,uz\}$ with $\{y_1z, y_2u\}$, we obtain a graph $G_9\in\langle D \rangle$, which contains the theta graph $\Theta(1,2,2)$ as subgraph.

Suppose $z$ is adjacent to $x_1$. By switching $\{x_1z, y_1y_2,y_2u\}$ with $\{y_1x_1, y_1z, y_1y_2\}$, we obtain a graph $G_{10}\in\langle D \rangle$, which contains the theta graph $\Theta(1,2,2)$ as subgraph.

Suppose $z$ is adjacent to $y_1$. By switching $\{y_1z,y_2v\}$ with $\{y_1v, y_2z\}$, we obtain a graph $G_{11}\in\langle D \rangle$, which contains the theta graph $\Theta(1,2,2)$ as subgraph.

If $t=4$, then let $uv$, $ux_1v$ and $uy_1y_2y_3v$ be the three paths of $\Theta(1,2,4)$. Since $n\geq7$, there is a vertex $z\in V(G)\setminus V(\Theta(1,2,4))$ such that $z$ is adjacent to one vertex in $\Theta(1,2,4)$. By symmetry, we only consider that  $z$ is adjacent to $u$, or $x_1$, or $y_1$, or $y_2$.  

Suppose $z$ is adjacent to $u$. By switching $\{y_2y_3,uz\}$ with $\{y_2z, y_3u\}$, we obtain a graph $G_{12}\in\langle D \rangle$, which contains the theta graph $\Theta(1,2,2)$ as subgraph.

Suppose $z$ is adjacent to $x_1$. By switching $\{x_1z, y_1y_2,y_3v\}$ with $\{y_1v, y_2x_1, y_3y_z\}$, we obtain a graph $G_{13}\in\langle D \rangle$, which contains the theta graph $\Theta(1,2,2)$ as subgraph.

Suppose $z$ is adjacent to $y_1$. By switching $\{y_1z,y_3v\}$ with $\{y_1v, y_3z\}$, we obtain a graph $G_{14}\in\langle D \rangle$, which contains the theta graph $\Theta(1,2,2)$ as subgraph.

Suppose $z$ is adjacent to $y_2$. By switching $\{y_1u, y_2y_3,y_2z,y_3v\}$ with $\{y_1y_3, y_2u, y_2v, y_3z\}$, we obtain a graph $G_{15}\in\langle D \rangle$, which contains the theta graph $\Theta(1,2,2)$ as subgraph.
$\hfill\Box$

\begin{lem}\label{lem49}
Let $D=(d_1,d_2,\ldots,d_n)$ be a forcibly bicyclic graphic sequence. If $n\geq7$, $D\neq(4,2^6)$ and $D\neq(5,2^6,1)$, then $d_7=1$.
 \end{lem}

\noindent{\bf Proof.}
By Lemma \ref{lem48}, there is a graph $G\in \langle D \rangle$ such that $G$ contains either the bowtie $B(3,3)$ or the kite $\Theta(1,2,2)$ as subgraph.  
On the contrary, assume $d_7\ge 2$. 

\noindent{\bf Case 1.}  $G$ contains the bowtie $B(3,3)$.

Let $ux_1x_2u$ and $uy_1y_2u$ be the two cycles of $B(3,3)$. Since $n\geq7$, there are two vertices  $z_1, z_2\in V(G)\setminus V(B(3,3))$ such that their degrees are all greater than 1. If $z_1$ and $z_2$ belong to $Tw_1-w_1$ and $Tw_2-w_2$, respectively, for two different  vertices $w_1, w_2\in V(B(3,3))$, then  there is one pendant vertex in  $Tw_1-w_1$ and one pendant vertex in  $Tw_2-w_2$ such that their distance is at least 5, which contradicts to Corollary \ref{cor21}.  If $z_1$ and $z_2$ belong to the same $Tw_3-w_3$ for some $w_3\in V(B(3,3))$ but they are in different components of $Tw_3-w_3$, then  there is one pendant vertex in the component  $Tw_1-w_1$ containing $z_1$ and one pendant vertex in  the component $Tw_2-w_2$ containing $z_2$ such that their distance is at least 4, which contradicts to Corollary \ref{cor21}.  So we assume  $z_1$ and $z_2$ belong to the same component of  $Tw-w$ for some $w\in V(B(3,3))$.  then there is a path of length three in $Tw$ from $w$, say $wv_1v_2v_3$. By symmetry, we only consider $w=x_1$ or $w=u$.

If $w=x_1$, then by switching $\{y_1u,y_2u, v_1v_2, v_2v_3\}$ with $\{y_1v_2,y_2v_2, v_1u,v_3u\}$, we obtain a disconnected realization of $D$, a contradiction.

If $w=u$, then by  $D\neq(5,2^6,1)$, there is a vertex $a\in V(G)\setminus (V(B(3,3))\cup\{v_1,v_2,v_3\})$ such that $a$ is adjacent to some vertex in $V(B(3,3))\cup\{v_1,v_2,v_3\}$. If $a$ is adjacent to one vertex of $V(B(3,3))$, then we can find two pendant vertex with distance at least 4, which contradicts to Corollary \ref{cor21}. If $a$ is adjacent to $v_1$, then by switching $\{x_1x_2, y_1y_2, av_1, v_1v_2\}$ with $\{x_2y_2, v_1x_1,v_1y_1,au_2\}$, we obtained a disconnected realization of $D$, a contradiction. If $a$ is adjacent to $v_2$, then by switching $\{x_1x_2, y_1y_2, av_2, v_2v_3\}$ with $\{x_2y_2, v_2x_1,v_2y_1,au_3\}$, we obtained a disconnected realization of $D$, a contradiction. If $a$ is adjacent to $v_3$, then by switching $\{uv_1, av_3\}$ with $\{au, v_1v_3\}$, we obtained a disconnected realization of $D$, also a contradiction. 

\noindent{\bf Case 2.}  $G$ contains  the theta graph $\Theta(1,2,2)$.

Let $uv$, $ux_1v$ and $uy_1v$ be the three paths of $\Theta(1,2,2)$. Since $n\geq7$, there are three vertices $z_1,z_2,z_3\in V(G)\setminus V(\Theta(1,2,2))$ such that their degrees are all greater than 1. 
If $z_1, z_2, z_3$ belong to $Tw-w$ for some $w\in V(\Theta(1,2,2))$, then  $Tw$ is not a star or double-star, which contradicts to Corollary \ref{cor22}. So there are two vertices in $V(\Theta(1,2,2))$, say $w_1$ and $w_2$, such that $V(Tw_1-w_1)\cap \{z_1, z_2, z_3\}\neq \emptyset$  and $V(Tw_2-w_2)\cap \{z_1, z_2, z_3\}\neq \emptyset$. Then there are one pendant vertex  in $Tw_1-w_1$ and one pendant vertex in  $Tw_2-w_2$ such that their distance is at least 5, which also contradicts to Corollary \ref{cor22}. 
$\hfill\Box$

\begin{lem}\label{lem410}
Let $D=(d_1,d_2,\ldots,d_n)$ be a forcibly bicyclic graphic sequence. If $n\geq7$, then $d_4=2$.
 \end{lem}

\noindent{\bf Proof.} If $D=(4,2^6)$, then the result holds. So assume $D\neq(4,2^6)$. 
By Lemma \ref{lem48}, there is a graph $G\in \langle D \rangle$ such that $G$ contains either the bowtie $B(3,3)$ or the kite $\Theta(1,2,2)$ as subgraph. So $d_4\geq2$.  
On the contrary, assume $d_4\ge 3$. 

\noindent{\bf Case 1.}  $G$ contains the bowtie $B(3,3)$.

Let $ux_1x_2u$ and $uy_1y_2u$ be the two cycles of $B(3,3)$. If there are a vertex in $\{x_1,x_2\}$ and a vertex $y\in \{y_1,y_2\}$ such that $d(x)\geq3$ and $d(y)\geq3$, then we can find a pendant vertex in  $Tx-x$ and a pendant vertex in  $Ty-y$ such that their distance is at least 4, which contradicts to Corollary \ref{cor21}. So $d(x_1)=d(x_2)=2$ or $d(y_1)=d(y_2)=2$. Without loss of generality, assume $d(x_1)=d(x_2)=2$. 
Since $d_4\ge 3$, there exists a vertex $z$ with degree at least three in $Ty_1-y_1$, or $Ty_2-y_2$, or $Tu-u$. 
If $z$ is in $Ty_1-y_1$, let $z_1$ and $z_2$ be two neighbors of $z$ not on the path between $y_1$ and $z$. Then by switching $\{x_1x_2,x_1u, zz_1,zz_2\}$ with $\{x_1z_1,x_1z_2,zx_2,zu\}$, we obtain a disconnected realization of $D$, a contradiction. Similarly, if $z$ is in $Ty_2-y_2$, then we can also obtain a disconnected realization of $D$, a contradiction. It is remaining to consider the case that $z$ is in $Tu-u$. Suppose that  $d(y_1)\geq3$ or $d(y_2)\geq3$.  Then we can find a pendant vertex in  $Tu-u$ and a pendant vertex in  $Ty_1-y_1$ or  $Ty_2-y_2$ such that their distance is at least 4, which contradicts to Corollary \ref{cor21}. So assume $d(y_1)=d(y_2)=2$. Then $Tu$  contains at least three vertices with degree at least three, which contradicts to Corollary \ref{cor22}.

\noindent{\bf Case 2.}  $G$ contains  the theta graph $\Theta(1,2,2)$.

Let $uv$, $ux_1v$ and $uy_1v$ be the three paths of $\Theta(1,2,2)$. If  $d(x_1)\geq3$ and $d(y_1)\geq3$, then we can find a pendant vertex in  $Tx_1-x_1$ and a pendant vertex in  $Ty_1-y_1$ such that their distance is at least 4, which contradicts to Corollary \ref{cor21}. So $d(x_1)=2$ or $d(y_1)=2$. Without loss of generality, assume $d(x_1)=2$. Since $d_4\ge 3$, there exists a vertex $z$ with degree at least three in $Ty_1-y_1$, or $Tu-u$, or $Tv-v$. 

If $z$ is in $Ty_1-y_1$, let $z_1$ and $z_2$ be two neighbors of $z$ not on the path between $y_1$ and $z$. Then by switching $\{x_1u,x_1v, zz_1,zz_2\}$ with $\{x_1z_1,x_1z_2,zu,vz\}$, we obtain a disconnected realization of $D$, a contradiction. 

If $z$ is in $Tu-u$, then $d(y_1)=2$. Otherwise, we can find a pendant vertex in $Tu-u$ and a pendant vertex in $Ty_1-y_1$ such that their distance is at least 4, which contradicts to Corollary \ref{cor21}. So there is another vertex $z'$ other than $z$ such that it has degree at least 3 and belongs to  $Tu-u$, or $Tv-v$.  If $z$ and $z'$ belong to different components of $Tu-u$ or $Tv-v$, then we can find two pendant vertices such that their distance is at least 4, which contradicts to Corollary \ref{cor21}.  So $z$ and $z'$ belong to one component of $Tu-u$. Then one of $z$ and $z'$,  say $z$, is not adjacent to $u$. Let $z_3$ and $z_4$ be two neighbors of $z$ not on the path between $u$ and $z$. Then by switching $\{x_1u,x_1v, zz_3,zz_4\}$ with $\{x_1z_3,x_1z_4,zu,vz\}$, we obtain a disconnected realization of $D$, a contradiction. 
The case $z$ is in $Tv-v$ can be proved similarly.
$\hfill\Box$

\begin{thm}
Let $D=(d_1,d_2,\ldots,d_n)$ be a non-negative and non-increasing integer sequence. Then
$D$ is a forcibly bicyclic graphic sequence if and only if one of the following holds:

$(1)$ $D=(3^2,2^4)$, or $D=(3^3,2^2,1)$, or $D=(3^3,2^3,1)$, or $D=(4,2^6)$, or $D=(4,3,2^4,1)$, or $D=(5,2^6,1)$;

$(2)$ $D=(n-1,2^4,1^{n-5})$,  where $n\ge5$; 

$(3)$ $D=(n-2,2^5,1^{n-6})$, where $n\ge6$;

$(4)$ $D=(n-2,3,2^3,1^{n-5})$, where $n\ge5$;

$(5)$ $D=(r,s,t,2,1^{n-4})$, where $n\ge4$, $r\ge s\ge 3$, $t\ge 2$ and $r+s+t=n+4$.
\end{thm}

\noindent{\bf Proof.}
First, we show the if part. Clearly, all of the sequences in (1)-(5) are graphic sequences and satisfy $\Sigma_{i=1}^n d_i=2n+2$. It is sufficient to prove that $G$ is connected for any graph $G\in \langle D \rangle$.

$(1)$ By Lemma \ref{lem47},  $D=(3^2,2^4)$ and $D=(3^3,2^2,1)$ are forcibly bicyclic graphic sequences.

Let $G\in \langle(3^3,2^3,1)\rangle$. Assume $d(u)=d(v)=d(w)=3$.  Since the degrees of vertices of $V(G)\setminus\{u,v,w\}$ are at most 2, $\{u,v,w\}$ are not an independent set of $G$. Without loss of generality, assume $uv\in E(G)$. 
If $|N(u)\cap N(v)|=0$, then there exists exactly one vertex $x$ in $V(G)\setminus (N[u]\cup N[v])$. By $d(x)\geq1$, $x$ must be adjacent to one vertex in $N(u)\cup N(v)$. Thus $G$ is connected. 
If $|N(u)\cap N(v)|=1$, there are exactly two vertices $x$ and $y$ in $V(G)\setminus (N[u]\cup N[v])$.  Since there is only one vertex with degree 1 in $G$, one of $x$ and $y$,  say $x$, has degree at least 2. Then $x$ is adjacent to some vertex in $N(u)\cap N(v)$ and $y$ is adjacent to some vertex in $N(u)\cap N(v)\cup\{x\}$.  Therefore, the graph $G$ is connected.
If $|N(u)\cap N(v)|=2$, there are exactly three vertices $x$, $y$ and $z$ in $V(G)\setminus (N_G[u]\cup N_G[v])$.  The degrees of $x$, $y$ and $z$ are at most two, otherwise there would be at least four vertices with degree 3. Assume, without loss of generality, $d(x)=d(y)=2$ and $d_G(z)=1$. Then one of $x$ and $y$, say $x$ is adjacent to one vertex in $N(u)\cap N(v)$, and $y$ is adjacent to some vertex in $N(u)\cap N(v)\cup\{x\}$, and $z$ is adjacent to some vertex in $N(u)\cap N(v)\cup\{x,y\}$. Therefore, the graph $G$ is connected.

Let $G\in \langle (4,2^6) \rangle$. Assume $d(u)=4$, $N(u)=\{u_1,u_2,u_3,u_4\}$ and $V(G)\setminus N[u]=\{x,y\}$.  Since $d(x)=d(y)=2$, both $x$ and $y$ is adjacent to some vertex in  $\{u_1,u_2,u_3,u_4\}$.
Therefore, the graph $G$ is connected.

Let $G\in \langle(4,3,2^4,1)\rangle$. Assume $d(u)=4$, $N(u)=\{u_1,u_2,u_3,u_4\}$ and $V(G)\setminus N[u]=\{x,y\}$. Since there is only one vertex with degree 1 in $G$, one of $x$ and $y$,  say $x$, has degree at least 2. Then $x$ is adjacent to some vertex in $N(u)\cap N(v)$ and $y$ is adjacent to some vertex in $N(u)\cap N(v)\cup\{x\}$.  Therefore, the graph $G$ is connected.

Let $G\in \langle (5,2^6,1)\rangle$. Assume $d(u)=5$ and  $V(G)\setminus N[u]=\{x,y\}$. Since there is only one vertex with degree 1 in $G$, one of $x$ and $y$,  say $x$, has degree at least 2. Then $x$ is adjacent to some vertex in $N(u)$ and $y$ is adjacent to some vertex in $N(u)\cup\{x\}$.  Therefore, the graph $G$ is connected.

$(2)$  Let $G\in\langle(n-1,2^4,1^{n-5})\rangle$. Since $G$ has one vertex with degree $n-1$, the graph $G$ is clearly connected.
 
(3)-(4) Let $G\in \langle (n-2,2^5,1^{n-6})\rangle$ or $G\in \langle(n-2,3,2^3,1^{n-5})\rangle$. Assume $d(u)=n-2$.
Then there exists exactly one vertex $x$ in $V(G)\setminus N[u]$. By $d(x)\geq1$, $x$ must be adjacent to one vertex in $N(u)$.  Therefore, the graph $G$ is connected.

$(5)$ Let $G\in \langle (r,s,t,2,1^{n-4})\rangle$ and $d(x)=r, d(y)=s$, $d(z)=t$, $d(w)=2$. Since all vertices in $V(G)\setminus\{x,y,z,w\}$ have degree 1, it follows that $x$, $y$ and $z$ must be adjacent to at least $r+s+t-8$ vertices with degree one. Since $r+s+t=n+4$ and $G$ has exactly $n$
vertices, we conclude that $x$, $y$ and $z$ are adjacent to exactly $n-4$ vertices with degree 1. Therefore, $G$ is the graph obtained from the kite $\Theta(1,2,2)$ with vertex set $\{x,y,z,w\}$ by adding $r-3$, $s-s$ and $t-2$ pendant vertices to $x$, $y$ and $z$, respectively, where $x$ and $y$ have degree three in the kite $\Theta(1,2,2)$. Thus,  $G$ is connected.

Next, we prove the only if part. Let $D$ be a forcibly bicyclic graphic sequence. Then $\Sigma_{i=1}^n d_i=2n+2$, $d_1\geq3$, and $d_2\geq d_3\geq d_4\geq2$.
If $n\leq6$, then, by Lemma \ref{lem47},  the result holds. In addition, $D=(4,2^6)$ and $D=(5,2^6,1)$ belong to (1). So we assume $n\geq7$, $D\neq(4,2^6)$ and $D\neq(5,2^6,1)$
in the following. By Lemmas \ref{lem49} and \ref{lem410}, we have $d_7=1$ and $d_4=2$.
Set $d_1=r$, $d_2=s$ and $d_3=t$.   We consider two cases as follows.

\noindent{\bf Case 1.} $d_6=1$.

By $d_4 = 2$,  we  have  $d_5=1$ or 2.

\noindent{\bf Subcase 2.1.}  $d_5=1$.

Since $\Sigma_{i=1}^n d_i=2n+2$, we have $r+s+t=2n+2-2-(n-4)=n+4$.  Thus $D=(r,s,t,2,1^{n-4})$, where $r+s+t=n+4$. 
If $s=2$, then $t=2$ and $D=(n,2^3,1^{n-4})$, which contradicts $d_1\leq n-1$. 
So $s\ge3$ and $D$ satisfies (5).

\noindent{\bf Subcase 2.2.} $d_5=2$.

Since $\Sigma_{i=1}^n d_i=2n+2$, we have $r+s+t=2n+2-2-2-(n-5)=n+3$.  Thus $D=(r,s,t,2^2,1^{n-5})$, where $r+s+t=n+3$. 

If $s=2$, then $t=2$ and $r=n-3$. So  $D=(n-1,2^4,1^{n-5})$, which satisfies (2). 

If $s=3$, then $D=(n-2,3,2^3,1^{n-5})$ or $D=(n-3,3^2,2^2,1^{n-5})$. Since $n\geq7$, the graph obtained from the disjoint union of $K_{1,n-3}$ and $K_2$ by adding three edges to $K_{1,n-3}$ (where those three edges induce a path on four vertices) has degree sequence $(n-3,3^2,2^2,1^{n-5})$, but is not connected, a contradiction. Thus $D=(n-2,3,2^3,1^{n-5})$, which belongs to (4).

Now we assume $s\geq4$. 
Then the graph obtained from  the disjoint union of $\Theta(2,2,2)+uv$ (where $u$ and $v$ are two vertices with degree 3 in $\Theta(2,2,2)$) and $K_2$ by adding $r-4$, $s-4$ and $t-2$ pendant vertices to the two vertices with degree 4 and one vertex with degree 2 in $\Theta(2,2,2)+uv$, respectively, has degree sequence $(r,s,t,2^2,1^{n-5})$. But it is not connected, which  contradicts that $D$ is a bicyclic graphic sequence.

\noindent{\bf Case 2.} $d_6=2$.

By $d_4=2$, we have   $D=(r,s,t,2^3,1^{n-6})$, where $r+s+t=n+2$. 

If $s=2$, then $t=2$ and $r=n-2$. So  $D=(n-2,2^5,1^{n-6})$, which satisfies (3). 

If $s=3$, then $D=(n-3,3,2^4,1^{n-6})$ or $D=(n-4,3^2,2^3,1^{n-6})$. The graph $G_1$ obtained from the disjoint union of $K_{1,n-3}$ and $K_2$ by adding three edges to $K_{1,n-3}$ (where those three edges induce the union of  a path on three vertices and a path on two vertices) has degree sequence $(n-3,3,2^4,1^{n-6})$ but is not connected, a contradiction. So $D\neq(n-3,3,2^4,1^{n-6})$. The graph $G_2$ obtained from the disjoint union of $K_{1,n-4}$ and $K_{1,2}$ by adding three edges to $K_{1,n-4}$ (where those three edges induce a path on four vertices) has degree sequence $(n-4,3^2,2^3,1^{n-6})$ but is not connected, a contradiction. So $D\neq (n-4,3^2,2^3,1^{n-6})$.

Now we assume $s\geq4$. 
Then the graph obtained from  the disjoint union of $\Theta(2,2,3)+uv$ (where $u$ and $v$ are two vertices with degree 3 in $\Theta(2,2,3)$) and $K_2$ by adding $r-4$, $s-4$ and $t-2$ pendant vertices to the two vertices with degree 4 and one vertex with degree 2 in $\Theta(2,2,3)+uv$, respectively, has degree sequence $(r,s,t,2^3,1^{n-6})$. But it is not connected, which contradicts that $D$ is a bicyclic graphic sequence.
$\hfill\Box$

\end{document}